\documentclass[11pt]{amsart}
\usepackage{epsfig}

\DeclareMathOperator{\conv}{conv}

\DeclareMathOperator{\spanOp}{span} \DeclareMathOperator{\adj}{adj}

\newtheorem{theorem}{Theorem}

\newtheorem{lemma}[theorem]{Lemma}

\newtheorem{definition}[theorem]{Definition}

\def\Prob{\mathrm{Prob}}
\def\rr{\mathbb{R}}

\providecommand{\abs}[1]{\left\lvert#1\right\rvert}

\providecommand{\Z}{\mathbb{Z}} \providecommand{\R}{\mathbb{R}}

\author{Sergi Elizalde}
\address{Department of Mathematics, Dartmouth College, Hanover, NH 03755}
\email{sergi.elizalde@dartmouth.edu}
\author{Kevin Woods}
\address{Department of Mathematics, University of California, Berkeley, CA 94720}
\email{kwoods@math.berkeley.edu}

\title{The Probability of Choosing Primitive Sets}

\begin{document}
\maketitle

\begin{abstract}
We generalize a theorem of Nymann that the density of points in
$\Z^d$ that are visible from the origin is $1/\zeta(d)$, where
$\zeta(a)$ is the Riemann zeta function $\sum_{i=1}^{\infty}1/i^a$.
A subset $S\subset \Z^d$ is called primitive if it is a $\Z$-basis
for the lattice $\Z^d\cap\spanOp_\R(S)$, or, equivalently, if $S$
can be completed to a $\Z$-basis of $\Z^d$.  We prove that if $m$
points in $\Z^d$ are chosen uniformly and independently at random
from a large box, then as the size of the box goes to infinity, the
probability that the points form a primitive set approaches
$1/[\zeta(d)\zeta(d-1)\cdots\zeta(d-m+1)]$.
 \end{abstract}

\section{Introduction}

A classic result in number theory is that, if a point in $\Z^2$ is
chosen ``at random,''  the probability that the point is visible
from the origin (that is, not hidden by another point in $\Z^2$) is
$\frac{1}{\zeta(2)}$, where $\zeta(a)$ is the Riemann zeta function
$\sum_{i=1}^{\infty}\frac{1}{i^a}$ (see \cite{Apo} for a proof using
Euler's totient function).  More precisely, for a given $n$, if we
choose an integer point $(a,b)$ uniformly at random from the box
$[-n,n]\times[-n,n]$ and compute the probability that $(a,b)$ is
visible from the origin, then as $n$ approaches infinity, this
probability approaches $\frac{1}{\zeta(2)}$.

J.E. Nymann generalized this result to higher dimensions \cite{Nym}:
if a point in $\Z^d$ is chosen at random, then the probability that
the point is visible from the origin is $\frac{1}{\zeta(d)}$.  This
theorem is true for $d\ge 2$ and is, in effect, true for $d=1$:  the
only points in $\Z^1$ that are visible from the origin are $\pm 1$,
so the probability is 0, and $\zeta(1)$ diverges so that
$\frac{1}{\zeta(1)}=0$.

An obvious way to restate the condition that a point
$s=(a_1,a_2,\ldots,a_d)\in\Z^d$ is visible from the origin is that
$\gcd(a_1,\ldots,a_d)=1$.  We will restate the condition in a
lattice theoretic context, so that it may be generalized to picking
more than one point in $\Z^d$. A point $s$ is visible from the
origin if and only if $\{s\}$ is a $\Z$-basis for the lattice
$\spanOp_{\rr}(s)\cap\Z^d$.  In general, given a set
$S=\{s_1,s_2,\ldots,s_m\}\subset\Z^d$, where $1\le m\le d$, we say
that $S$ is \emph{primitive} if $S$ is a $\Z$-basis for the lattice
$\spanOp_{\rr}(S)\cap\Z^d$.  An equivalent definition \cite{Lek} is
that $S$ is primitive if and only if $S$ can be completed to a
$\Z$-basis of all of~$\Z^d$.

In this paper we prove that if $S$ is chosen ``at random,'' then the probability that $S$ is primitive is
\[\frac{1}{\zeta(d)\zeta(d-1)\cdots\zeta(d-m+1)}.\]
To be precise, we prove the following theorem.

\begin{theorem}
\label{Theorem:Main} Let $d$ and $m$ be given, with $m<d$.  For
$n\in\Z_+$, $1\le k\le m$, and $1\le i\le d$, let $b_{n,k,i}\in\Z$.
For a given $n$, choose integers $s_{ki}$ uniformly (and
independently) at random from the set $b_{n,k,i}\le s_{ki}<
b_{n,k,i}+n$.  Let $s_k=(s_{k1},\ldots,s_{kd})$ and let
$S=\{s_1,s_2,\ldots,s_m\}$.

If  $\abs{b_{n,k,i}}$ is bounded by a polynomial in $n$, then, as $n$ approaches infinity, the
probability that $S$ is a primitive set approaches
\[\frac{1}{\zeta(d)\zeta(d-1)\cdots\zeta(d-m+1)},\]
where $\zeta(a)$ is the Riemann zeta function
$\sum_{i=1}^{\infty}\frac{1}{i^a}$.
\end{theorem}

When $m=1$, this theorem gives the classic result ($d=2$) and
Nymann's result. Note also that, if $m=d$ and we choose $S$ of size
$m$, then the probability that $S$ is primitive (i.e., that it is a
basis for $\Z^d$) approaches zero. This agrees with the theorem in
the sense that we would expect the probability to be
\[\frac{1}{\zeta(d)\zeta(d-1)\cdots\zeta(1)},\]
but $\zeta(1)$ does not converge.

The statement of the theorem uses more general boxes than $[-n,n]^d$
to pick the $s_k$ from.  We do this because the more general result
is needed in \cite{EW1}.  That paper was the original inspiration
for this theorem: we discovered it in an attempt to prove a fact in
computational biology and Bayesian network theory.  Since the
concept of primitive sets is important in the geometry of numbers,
we are proving this theorem in this separate paper.  Note that some
bound on the $b_{n,k,i}$ in terms of $n$ is needed; otherwise one
could use the Chinese Remainder Theorem to construct arbitrarily
large boxes from which \emph{no} primitive sets could be selected
(even for $d=2$, $m=1$).

In Section 2, we present an outline of the proof.  The outline is a
full proof in every respect, except that we ignore the error
estimations in our probabilities.  In that sense, it is the
``moral'' proof of the result.  In Section 3, we fill the holes by
proving that the error estimates approach zero as $n$ approaches
infinity.  The methods in Section 3 are themselves of interest,
using concepts from triangulations of point sets, the metric
geometry of polytopes (cross-sections of d-cubes), analytic number
theory (consequences of the Prime Number Theorem), and the geometry
of numbers.

\section{Outline of the proof}

We proceed by induction on $m$.

If $m=0$, the theorem is trivially true. Assume that
the theorem is true for $m-1$, and we will prove it for $m$.  The
probability that $S=\{s_1,s_2,\ldots,s_m\}$ is primitive is the product
\begin{align*}
\Prob&_{\mathcal{P}_n}\big(\{s_1,\ldots,s_{m-1}\} \text{ is
primitive}\big)\\
&\cdot\Prob_{\mathcal{P}_n}\big(S \text{ is primitive, given
that }\{s_1,\ldots,s_{m-1}\} \text{ is primitive}\big),
\end{align*}
where
$\mathcal{P}_n$ is the probability distribution, for a given $n$,
from which we are choosing $S$. The first term in the product
approaches
\[\frac{1}{\zeta(d)\zeta(d-1)\cdots\zeta(d-m+2)},\]
as $n\rightarrow\infty$, by the inductive hypothesis, so we must
show that the second term approaches $\frac{1}{\zeta(d-m+1)}$.

Indeed, suppose $\{s_1,\ldots,s_{m-1}\}$ is given and is primitive,
and we choose $s_m=(s_{m1},\ldots,s_{md})$ (independently from the
other $s_i$) according to the probability distribution
$\mathcal{P}_n$. Let $A$ be the $(m-1)\times d$ integer matrix whose
rows are $s_1,\ldots,s_{m-1}$.  We will need the following lemma, to
find a simpler matrix whose rows also form a primitive set.

\begin{lemma}
\label{lemma:U} Let $A$ be a matrix in $\Z^{p\times q}$, and let $U$
be a unimodular matrix (i.e., $\det(U)=\pm 1$) in $\Z^{q\times q}$.
The rows of $A$ form a primitive set if and only if the rows of $AU$
also form a primitive set.
\end{lemma}

\begin{proof}
Suppose the rows of $A$ form a primitive set.  Let $a\in\Z^q$ be in
the $\R$-span of the rows of $AU$, that is, $a=xAU$, where $x$ is a
matrix in $\R^{1\times p}$.  In order to show that the rows of $AU$
form a primitive set, we must show that $x$ is actually integral.
Indeed, $aU^{-1}=xA\in\Z^q$ is in the $\R$-span of the rows of $A$,
and since the rows of $A$ form a primitive set, $x$ must integral. This also proves the converse, as $U^{-1}$ is unimodular and
$A=(AU)U^{-1}$.
\end{proof}

The matrix $U$ we will choose is a matrix that puts $AU$ into \emph{Hermite normal form}.

\begin{definition}
A matrix $B\in\Z^{p\times q}$ is in \emph{Hermite normal form} if
\begin{enumerate}
\item $B_{ij}=0$ for all $j>i$,
\item $B_{ii}> 0$ for all $i$, and
\item $0\le B_{ij} < B_{ii}$ for all $j<i$.
\end{enumerate}
\end{definition}

Given any integer matrix $B$ of full row rank, there exists a
unimodular matrix $U$ such that $BU$ is in Hermite normal form (see,
e.g., \cite{GLS}; $U$ will not, in general, be unique).  This fact,
together with the following lemma, gives a convenient
characterization of when $S$ is a primitive set.

\begin{lemma}
\label{Lemma:RelPrime}
Let $\{s_1,\ldots,s_{m-1}\}\subset\Z^d$ be a
primitive set, and let $s_m\in\Z^d$ be given. Let $A$ be the (full
row rank) matrix with rows $s_1,\ldots,s_{m-1}$, and let $U$ be a
matrix such that $AU$ is in Hermite normal form.  Let $U^{(i)}$ be
the $i$-th column of $U$. Then $\{s_1,\ldots,s_m\}$ is a primitive
set if and only if the $s_mU^{(i)}$, for $m\le i\le d$, are
relatively prime.
\end{lemma}

\begin{proof}
By Lemma \ref{lemma:U}, the rows of $AU$ form a primitive set. It
follows that $(AU)_{ii}=1$, for $1\le i\le m-1$ (otherwise $e_i$,
the $i$-th standard basis vector, would be in the $\R$-span of the
rows of $AU$, but not in the $\Z$-span). Then, from the definition
of Hermite normal form, $(AU)_{ij}=0$ for $i\ne j$.  Let $A'$ be the
matrix with rows $s_1,\ldots,s_{m}$ (that is, $A'$ is $A$ with the
additional row $s_m$ appended).  By Lemma \ref{lemma:U},
$\{s_1,\ldots,s_m\}$ is a primitive set if and only if the rows of
$A'U$ form a primitive set.  We see that this is true if and only if
the $(A'U)_{mi}$, for $m\le i\le d$, are relatively prime (indeed,
the index of the lattice $\spanOp_\Z\{s_1,\ldots,s_m\}$ within
$\Z^{d}\cap\spanOp_{\R}\{s_1,\ldots,s_m\}$ is $\gcd\{(A'U)_{mi}:\
m\le i\le d\}$). Since $(A'U)_{mi}=s_mU^{(i)}$, the lemma follows.
\end{proof}

Let $\mu:\Z_+\rightarrow \{-1,0,1\}$ be the M\"obius function
defined to be
\[\mu(D)=\begin{cases}
(-1)^i & \text{if $D$ is the product of $i$ distinct primes,}\\
0 & \text{if $D$ is divisible by the square of a prime.}
\end{cases}\]
Given $D\in\Z_+$, let $p_{nD}$ be the probability that $D$ divides
$s_mU^{(i)}$ for all $m\le i\le d$. Note that $p_{nD}$ is
independent of our choice of $U$, because, as we noted in the proof
of Lemma \ref{Lemma:RelPrime}, $\gcd\{s_mU^{(i)}:\ m\le i\le d\}$ is
the index of the lattice $\spanOp_\Z\{s_1,\ldots,s_m\}$ within
$\Z^{d}\cap\spanOp_{\R}\{s_1,\ldots,s_m\}$, which is independent of
$U$. Then, using inclusion-exclusion, the probability that the
$s_mU^{(i)}$, for $m\le i\le d$, are relatively prime is
\[\sum_{D=1}^{\infty}\mu(D)p_{nD}.\]

We expect each $p_{nD}$ to be approximately
$D^{-(d-m+1)}$.  In Section 3, we will show that
\begin{equation}
\label{Eqn:lim} \lim_{n\rightarrow
\infty}\sum_{D=1}^{\infty}\mu(D)p_{nD} =
\sum_{D=1}^{\infty}\mu(D)D^{-(d-m+1)}.
\end{equation}

Given that we have verified (\ref{Eqn:lim}), the following lemma
(applied to $a=d-m+1$) finishes the proof of the theorem.

\begin{lemma}
For any integer $a\ge 2$,
\[\sum_{D=1}^{\infty}\mu(D)D^{-a}=\frac{1}{\zeta(a)}.\]
\end{lemma}

\begin{proof}
Since $a\ge 2$, the sum is absolutely convergent, and we have that
\begin{align*}
\sum_{D=1}^{\infty}\mu(D)D^{-a}
&=\prod_{p \text{ prime}}\left(1-p^{-a}\right)\\
&=\frac{1}{\displaystyle\prod_{p \text{ prime}}\frac{1}{1-p^{-a}}}\\
&=\frac{1}{\displaystyle\prod_{p \text{ prime}}\left(1+p^{-a}+p^{-2a}+\cdots\right)}\\
&=\frac{1}{\displaystyle\sum_{i=1}^{\infty}i^{-a}}\\
&=\frac{1}{\zeta(a)}.
\end{align*}
\end{proof}

\section{Error Estimates}

The remaining piece of the proof is to demonstrate Equation
(\ref{Eqn:lim}), that is, that
\[\abs{\sum_{D=1}^{\infty}\mu(D)p_{nD} - \sum_{D=1}^{\infty}\mu(D)D^{-(d-m+1)}}\rightarrow 0\]
as $n\rightarrow \infty$.

We will need a bound on the entries of $U$, which the following
lemma will help us get.

\begin{lemma}
\label{lemma:HNFbound} Given a rank $p$ matrix $A\in\Z^{p\times q}$
and a bound $M_0$ such that $\abs{A_{ij}}< M_0$ for all $i,j$, there
exists a unimodular matrix $U$ such that
\begin{enumerate}
\item $AU$ is in Hermite normal form and
\item $\abs{U_{ij}}\le p!qM_0^p$ for all $i,j$.
\end{enumerate}
\end{lemma}

\begin{proof}
Let $B$ be the $q\times q$ matrix obtained be appending to $A$ the
rows $e_1, e_2,\ldots, e_{q-p}$ (where $e_i$ is the $i$-th standard
basis vector). Without loss of generality, we can assume that $B$ is
a nonsingular matrix (otherwise, we could have appended different
$e_i$).  Let $U$ be a unimodular matrix such that $BU$ is in Hermite
normal form. Note that $AU$ is also in Hermite normal form.

We will use the fact that
\begin{equation}\label{eqn:forU}
U=B^{-1}(BU)=\frac{1}{\det(B)}\adj(B)(BU), \end{equation} where
$\adj(B)$ is the adjugate (classical adjoint) of $B$, in order to bound the
entries of $U$. Since $BU$ is lower diagonal,
\[\abs{\det(B)}=\det(BU)=\prod_{i=1}^q(BU)_{ii}.\]
Therefore $(BU)_{ii}\le \abs{\det(B)}$ for all $i$, and, by the
definition of Hermite normal form, we conclude that $(BU)_{ij}\le
\abs{\det(B)}$ for all $i,j$.

Since the first $p$ rows of $B$ have entries bounded by $M_0$ and
the remaining rows are standard basis vectors, the entries of
$\adj(B)$ are bounded by $p!M_0^p$. Combining these two bounds, we
see that the entries of $\adj(B)(BU)$ are bounded by $q\cdot
p!M_0^p\cdot\abs{\det(B)}$.  Using (\ref{eqn:forU}) we conclude that
\[\abs{U_{ij}}\le \frac{1}{\abs{\det(B)}}q\cdot p!M_0^p\cdot\abs{\det(B)}=p!qM_0^p\]
for all $i,j$, as desired.
\end{proof}

Since the absolute value of the entries of $A$ are bounded by the
$b_{n,k,i}+n$, which we assume to be bounded by a polynomial in $n$,
Lemma \ref{lemma:HNFbound} shows that the unimodular matrix $U$ can
be chosen such that the absolute value of each entry of $U$ is
bounded by a polynomial in $n$.  This in turn implies that
$\abs{s_mU^{(i)}}$ is also bounded by a polynomial in $n$ (where
$U^{(i)}$ is the $i$-th column of $U$). Let $M=M(n)$ be our bound on
$\abs{s_mU^{(i)}}$; say $M$ is $O(n^k)$ for some $k$. Clearly, for
$D>M$, $p_{nD}=0$.

We have that

\begin{equation}
\label{eqn:error}
\begin{split}
&\abs{\sum_{D=1}^{\infty}\mu(D)p_{nD} - \sum_{D=1}^{\infty}\mu(D)D^{-(d-m+1)}}\\
&\ \ \le \abs{\sum_{D=1}^n\mu(D)\left(p_{nD}-D^{-(d-m+1)}\right)}
+\abs{\sum_{D=n+1}^M\mu(D)p_{nD}}\\
&\ \ \ \ \ \ +\abs{\sum_{D=M+1}^{\infty}\mu(D)p_{nD}}
+\abs{\sum_{D=n+1}^{\infty}\mu(D)D^{-(d-m+1)}}\\
&\ \ \le \sum_{D=1}^n\abs{p_{nD}-D^{-(d-m+1)}}
+\sum_{D=n+1}^Mp_{nD}
+0
+\sum_{D=n+1}^{\infty}D^{-(d-m+1)}.
\end{split}
\end{equation}

Of the three nonzero terms in the last expression,
$\sum_{D=n+1}^{\infty}D^{-(d-m+1)}$ certainly converges to zero as
$n$ approaches infinity, so it suffices to show that the first two
terms, $\sum_{D=1}^n\abs{p_{nD}-D^{-(d-m+1)}}$ and
$\sum_{D=n+1}^Mp_{nD}$, do as well.  We break our error computation
into these two cases.

Before we handle the two error sums in Lemmas \ref{lemma:1} and
\ref{lemma:2}, we set some common terminology.  Let $\mathcal{B}_n$
be the $d$-dimensional box of integers $\{s_m\in\Z^d:\ b_{n,m,i}\le
s_{mi} < b_{n,m,i} +n,\text{ for all }i\}$, which is the box from
which $s_m$ is chosen with uniform probability.  Given $D\in\Z_+$,
let $\Lambda_D\subset\Z^d$ be the lattice of integer vectors
$x\in\Z^d$ such that $D$ divides $x\cdot U^{(i)}$, for $m\le i\le
d$. $\Lambda_D$ is a sublattice of $\Z^d$ of index $D^{d-m+1}$. Let
$S_{nD}=\mathcal{B}_n\cap\Lambda_D$.  Then
\begin{equation}\label{eqn:p_nd}
p_{nD}=\frac{\abs{S_{nD}}}{n^d}.\end{equation}

\begin{lemma}\label{lemma:1}
As defined above,
\[\sum_{D=1}^n\abs{p_{nD}-D^{-(d-m+1)}}\]
converges to zero as $n\rightarrow \infty$.
\end{lemma}

\begin{proof}
Suppose $1\le D\le n$.  Let $L_D\subset\Z^d$ be the lattice of
integer vectors $(x_1,\ldots,x_d)\in\Z^d$ such that $D$ divides each
$x_i$.  $L_D$ is a sublattice of $\Z^d$ of index $D^d$.  In fact, we
see that $L_D$ is a sublattice of $\Lambda_D$, and therefore its
index in $\Lambda_D$ is $D^d/D^{d-m+1}=D^{m-1}$.

This means that if we look at any $D\times\cdots\times D$ cube, $C=\{(x_1,\ldots,x_d)\in\Z^d:\ r_i\le x_i < r_i+D\}$ for some $r_i\in\Z$ (that is, a translate of a fundamental parallelepiped of $L_D$), then $C$ contains exactly $D^{m-1}$ elements of $\Lambda_D$.  Since $\mathcal{B}_n$ can be
covered by $(\frac{n}{D}+1)^d$ such boxes, we have
that $\abs{S_D}\le D^{m-1}(\frac{n}{D}+1)^d$, and so
\[p_{nd}\le \frac{D^{m-1}(\frac{n}{D}+1)^d}{n^d}=\sum_{k=0}^d\binom{d}{k}\frac{D^{m-1-k}}{n^{d-k}}.\]
Similarly, $(\frac{n}{D}-1)^d$ disjoint $D\times\cdots\times D$ cubes can
be placed inside $\mathcal{B}_n$, and so
\[p_{nd}\ge \frac{D^{m-1}(\frac{n}{D}-1)^d}{n^d}=\sum_{k=0}^d\binom{d}{k}(-1)^{d-k}\frac{D^{m-1-k}}{n^{d-k}}.\]

Combining these two inequalities and moving the $k=d$ summand to the left-hand side, we see that
\[\abs{p_{nd}-\frac{1}{D^{d-m+1}}}\le
\sum_{k=0}^{d-1}\binom{d}{k}\frac{D^{m-1-k}}{n^{d-k}}\]
and so
\[
\sum_{D=1}^n\abs{p_{nD}-D^{-(d-m+1)}}\le\sum_{k=0}^{d-1}\left[\binom{d}{k}n^{-(d-k)}\sum_{D=1}^nD^{m-1-k}\right]
\]
which converges to zero as $n\rightarrow\infty$, proving the lemma.
\end{proof}

\begin{lemma}
\label{lemma:2}
As defined above,
\[\sum_{D=n+1}^M p_{nD},\]
 converges to zero as $n\rightarrow\infty$.
\end{lemma}

\begin{proof}
 Let
\[T_n=\bigcup_{D=n+1}^{M}S_{nD}.\]
Let $N_n$ be the maximum, over all $s_m\in\mathcal{B}_n$, of
\[\#\{D:\ n<D\le M \text{ and }s_m\in S_{nD}\}.\]
Then

\begin{align*}
\sum_{D=n+1}^M p_{nD}&=n^{-d}\sum_{D=n+1}^M\abs{S_{nD}}\\
&\le n^{-d}\abs{T_n}\cdot N_n
\end{align*}

We need to approximate $N_n$ and $\abs{T_n}$. We will repeatedly use
the following fact (see \cite{Apo}, p:294), which can be derived
from the Prime Number Theorem: for any $\epsilon>0$ and for any
$r\le M$, the number of factors of $r$ is $O(n^\epsilon)$ (more
precisely, for any $\delta>0$ and sufficiently large $r$, the number
of factors of $r$ is less than $r^{(1+\delta)\log2/\log\log r}$; now
we use that $r\le M$ is $O(n^k)$ for some $k$).

\ \\

\noindent \emph{Claim 1: $N_n$ is $O(n^{\epsilon})$.}\\

This follows immediately, as any element of the set
\[\{D:\ n<D\le M \text{ and }s_m\in S_{nD}\}\]
must be a factor of, say, $s_mU^{(m)}$, and this number has $O(n^{\epsilon})$ factors.

\ \\

\noindent \emph{Claim 2: $\abs{T_n}$ is $O(n^{d-\frac{1}{2}+\epsilon})$.}\\

Let $a=\gcd(U^{(i)}_1:\ m\le i\le d)$, where
$U^{(m)},U^{(m+1)},\ldots,U^{(d)}$ are the last $d-m+1$ columns of
$U$. Let $R$ be the set of integers greater than $n$ that are
factors of at least one of $a, 2a, 3a, \ldots,
\lfloor\sqrt{n}\rfloor a$. Each of the $\lfloor\sqrt{n}\rfloor$
numbers $i\cdot a$ such that $1\le i\le \lfloor\sqrt{n}\rfloor$ has
$O(n^\epsilon)$ factors, so $\abs{R}$ is
$O(n^{\frac{1}{2}+\epsilon})$.

We divide $T_n$ into two parts.  Let
\[T_{n1}=\bigcup_{D\in R}S_{nD}\]
and let $T_{n2}=T_n\setminus T_{n1}$.  We will show that both $\abs{T_{n1}}$ and $\abs{T_{n2}}$ are $O(n^{d-\frac{1}{2}+\epsilon})$, and so it will follow that $\abs{T_n}=\abs{T_{n1}}+\abs{T_{n2}}$ is also $O(n^{d-\frac{1}{2}+\epsilon})$.

\ \\

\noindent \emph{Claim 2a: $\abs{T_{n1}}$ is $O(n^{d-\frac{1}{2}+\epsilon}).$}\\

Given a $D\in R$, we want to estimate how large $S_{nD}$ is.
 Suppose first that $\conv(S_{nD})$ is a full dimensional polytope in
$\Z^d$, that is, its affine hull is all of $\R^d$. Triangulate
$\conv(S_{nD})$ into at least $\abs{S_{nD}}-d$ simplices whose
vertices are in $S_{nD}$ (this can always be done, see for example
\cite{DRS}).  Each simplex in the triangulation has volume at least
$\frac{1}{d!}D^{d-m+1}$, because the lattice $\Lambda_n$ (which
includes every point in $S_{nD}$) has index $D^{d-m+1}$ in $\Z^d$.
But $\conv(S_{nD})$ has volume at most $n^d$, because it lies in
$\mathcal{B}_n$.  Putting this together,
\[\frac{1}{d!}D^{d-m+1}(\abs{S_{nD}}-d)\le n^d,\]
and so
\[\abs{S_{nD}}\le d+d!\frac{n^d}{D^{d-m+1}}\le d+d!n^{m-1},\]
which is $O(n^{m-1})$.

On the other hand, if $\conv(S_{nD})$ is not full dimensional, then
let $k\le d-1$ be its dimension, and let $H$ be the $k$-dimensional
affine space such that $S_{nD}\subset H$.  The $k$-dimensional
Euclidean volume of $H\cap\mathcal{B}_n$ is at most
$\sqrt{2}^{d-k}n^{k}$, as proved in \cite{Bal}.  Again we can
triangulate $S_{nD}$ into at least $\abs{S_{nD}}-k$ simplices that
are $k$-dimensional. The best we can know this time is that each
simplex has volume at least $\frac{1}{k!}$.  Putting this together,
\[\frac{1}{k!}(\abs{S_{nD}}-k)\le \sqrt{2}^{d-k}n^{k},\]
and so $\abs{S_{nD}}$ is $O(n^k)$.

In either case, $\abs{S_{nD}}$ is $O(n^{d-1})$, and since $\abs{R}$
is $O(n^{\frac{1}{2}+\epsilon})$, $\abs{T_{n1}}$ is $O(n^{d-1}\cdot
n^{\frac{1}{2}+\epsilon})=O(n^{d-\frac{1}{2}+\epsilon})$.

\ \\

\noindent \emph{Claim 2b: $\abs{T_{n2}}$ is $O(n^{d-\frac{1}{2}+\epsilon})$.}\\

Recall that $a=\gcd(U^{(i)}_1:\ m\le i\le d)$. Without loss of
generality, we may assume that $U^{(m)}_1=a$ and $U^{(i)}_1=0$, for
$m+1\le i\le d$ (if not, we may perform elementary column operations
on the last $d-m+1$ columns of $U$ in order to put them in that
form; the matrix $AU$ will remain in Hermite normal form, because
the last $d-m+1$ columns of $AU$ are all zeros).  Note that $a<M$.

Now suppose $s_{m2},s_{m3},\ldots,s_{md}$ are given, such that
$b_{n,m,i}\le s_{mi}< b_{n,m,i}+n$.  Given $j$ such that
$b_{n,m,1}\le j<b_{n,m,1}+n$, define
\[t^{(j)}=(j,s_{m2},s_{m3},\ldots,s_{md}).\]
We will show that $O(n^{\frac{1}{2}+\epsilon})$ of the $t^{(j)}$ are
in $T_{n2}$ (for given $s_{m2},\ldots,s_{md}$).

Since $U^{(m+1)}_1=0$, $s':=t^{(j)}U^{(m+1)}$ is independent of $j$.
If $t^{(j)}\in S_{nD}$ for a particular $D$, then $D$ must be a
factor of $s'$, which has $O(n^{\epsilon})$ factors.  Therefore
there are only $O(n^{\epsilon})$ possible $D$ for which any of the
$t^{(j)}$ could be a member of $S_{nD}$.

Now let us consider, for a given $D\notin R$, how many of the
$t^{(j)}$ could be in $S_{nD}$.  If $t^{(j)}$ and $t^{(k)}$ are in
$S_{nD}$, then $D$ divides $t^{(j)}U^{(m)}$ and $t^{(k)}U^{(m)}$.
Therefore $D$ divides the difference
$t^{(j)}U^{(m)}-t^{(k)}U^{(m)}$, which is $(j-k)\cdot a$, since
$U^{(m)}_1=a$.  Since $D\notin R$, $D$ does not divide $a, 2a,
\ldots,\lfloor\sqrt{n}\rfloor a$, and so $\abs{j-k}>\sqrt{n}$.
Therefore the number of $j$ such that $t^{(j)}\in S_{nD}$ is at most
$n/\sqrt{n}=\sqrt{n}$.

Since there are $O(n^{\epsilon})$ possibilities for $D$, and since,
for a given $D\notin R$, the number of $t^{(j)}$ in $S_{nD}$ is
$O(n^{\frac{1}{2}})$, we conclude that $O(n^{\frac{1}{2}+\epsilon})$
of the $t^{(j)}$ are in $T_{n2}$.

Since there are $n^{d-1}$ choices for $s_{m2},\ldots,s_{md}$, we
have that $\abs{T_{n2}}$ is
\[O(n^{d-1}n^{\frac{1}{2}+\epsilon})=O(n^{d-\frac{1}{2}+\epsilon}).\]

\ \\

Combining our estimates of $N_n$ and $\abs{T_n}$, we have that
\begin{align*}
\sum_{D=n+1}^M p_{nD}&\le n^{-d}\abs{T_n}\cdot N_n\\
&=n^{-d}O(n^{d-\frac{1}{2}+\epsilon})O(n^{\epsilon})\\
&=O(n^{-\frac{1}{2}+2\epsilon}),
\end{align*}
and therefore $\sum_{D=n+1}^M p_{nD}$ converges to zero as $n$ approaches infinity.
\end{proof}

Combining Lemmas \ref{lemma:1} and \ref{lemma:2} with Equation
(\ref{eqn:error}), we have shown that
\[\abs{\sum_{D=1}^{\infty}\mu(D)p_{nD} - \sum_{D=1}^{\infty}D^{-(d-m+1)}}\rightarrow 0\]
as $n\rightarrow \infty$. This completes our error analysis and, together with Section 2, provides a complete proof of Theorem \ref{Theorem:Main}.

\end{document}